\pdfoutput=1
\documentclass{article}
\usepackage{header}
\usepackage{geometry}
\usepackage[backend=bibtex,style=numeric,doi=false,isbn=false,url=false,maxnames=10]{biblatex}
\DeclareMathOperator{\sgn}{sgn}

\AtBeginDocument{
}

\AtEndDocument{
  \printbibliography
  \bigskip
  \noindent\!\!\!\!\begin{tabular}{ll}
    Xiudi Tang & \quad \texttt{xiudi.tang@bit.edu.cn} \\
    Chengle Peng & \quad \texttt{chengle.peng@bit.edu.cn} \\
  \end{tabular} \smallskip \\
  {\footnotesize \textsc{\bitsms}}
}

\begin{document}

\title{Action-angle coordinates of spherical pendulums with symmetric quadratic potentials}
\author{Chengle Peng and Xiudi Tang}
\maketitle

\begin{abstract}
  We study the spherical pendulum system with an arbitrary potential function $V = V (z)$, which is an integrable system with a first integral whose Hamiltonian flow is periodic.
  We give an explicit solution to this integrable system and then we compute its action-angle coordinates.
  In the special case where the potential function is symmetric quadratic like $V = z^2$, we represent its action-angle coordinates in terms of elliptic integrals, and calculate the monodromy.
\end{abstract}

\section{Introduction} \label{sec:intro}

In Hamiltonian dynamics, there is a large class of dynamical systems called \emph{Liouville integrable systems} or just integrable systems for simplicity.
They have a maximal number of functionally independent quantities that are conserved under the dynamics.
Within the class of integrable systems, there is a special subclass whose elements display a rotational symmetry on one of the two components and are thus called \emph{semitoric systems}.

A semitoric system on a $4$-dimensional symplectic manifold $(M, \om)$ is a pair of Poisson commutative functions $(J, H)$ on $M$, where $J$ is the first integral generating an effective $S^1$-action, and $F = (J, H) \colon M \to \R^2$, viewed as a singular Lagrangian fibration, has singularities of a certain Morse--Bott type, with compact and connected fibers.
Four-dimensional integrable systems can have critical points of $6$ types: elliptic-regular, hyperbolic-regular, elliptic-elliptic, elliptic-hyperbolic, hyperbolic-hyperbolic, and focus-focus.
In the semitoric case, three of them remain: elliptic-regular, elliptic-elliptic, and focus-focus.
\'Alvaro Pelayo and San V\~u Ng\d{o}c in their papers \cite{MR2784664,MR2534101} classified simple semitoric systems up to fiber-preserving symplectomorphisms, where simplicity means that $J$-fiber contains at most one \emph{focus-focus singularity}. 
Removing the simplicity assumption, Joseph Palmer, \'Alvaro Pelayo, and the second named author in their paper \cite{MR4797864} generalized the classification theorem to all semitoric integrable systems.

The \emph{spherical pendulum} is one of the few natural Hamiltonian systems that is integrable and possesses focus-focus singularities; see \cref{eq:JH-in-cart-spdlm-gen}.
In a linear potential field (such as the gravity on the earth surface), there is a unique focus-focus point.
This is the very example for which Duistermaat \cite{MR596430} found that the \emph{action-angle coordinates} do not exist globally.
As a natural Hamiltonian system, its phase space being the cotangent bundle ensures the ease of quantization, and the kinetic energy being quadratic brings the elegance of the Legendre transformation to generate the Lagrangian function.
When we allow for a more general potential function, there can be multiple focus-focus points and other singularity behaviors.
In the case of $V (z) = z^2$ as in \cref{eq:JH-in-cart-spdlm-quad}, there are exactly two focus-focus points on the same fiber of the momentum map.
Efstathiou's book \cite{MR2152605} investigated this system, mostly on their Dirac--Poisson structures and topological properties.

In this paper, we take a step forward and compute a version of spherical pendulum whose potential is a general function $V (z)$ of the height $z$, solve the Hamiltonian flows of $J$ and $H$ explicitly in \cref{eq:sol-spdlm-gen}, and then write the action coordinates as a definite integral \cref{eq:action-spdlm-gen}.
In particular, we study thoroughly the case when $V (z) = z^2$.
We give explicit formulas for the joint flow by the action of the momentum map $F$ in \cref{eq:sol-spdlm-quad} and then for the action-angle coordinates \cref{eq:action-angle-sq-spdlm}, both in terms of elliptic integrals.
The monodromy \cref{eq:monodromy-sq-spdlm} we calculate is in accordance with the theory of focus-focus fibers in integrable systems.

There are only a few explicit semitoric systems so far whose symplectic invariants are actually computed.
Although the spherical pendulum is not strictly a semitoric system as $J$ is not proper, it falls into the class of faithful semitoric systems as in Hohloch--Sabatini--Sepe--Symington \cite{MR3843843}, and have similar semiglobal invariants as semitoric ones.
The spherical pendulum with a linear potential has a pinched torus fiber, while the one with a symmetric quadratic potential has a double-pinched torus fiber.
In particular, Dullin \cite{MR3017036} computed the Taylor series invariants of the spherical pendulum.
Other examples include Le Floch--Pelayo \cite{MR3927110} of coupled angular momenta and Alonso--Dullin--Hohloch \cite{MR3924558} of coupled spin-oscillators, but we are not aware of the computation for general spherical pendulums.
In the future, we plan to compute the Taylor series and the other invariants of general spherical pendulums so as to deepen the analysis on them. 

\section{Local coordinates on $T S^2$} \label{sec:local-coord}

It is helpful to embed the cotangent bundle $T^* S^2$ of the two-sphere into the cotangent bundle $T^* \R^3$ of the three-dimensional Euclidean space, with position coordinates $x, y, z$ and momentum coordinates $u, v, w$, equipped with the canonical symplectic form
\begin{equation*}
  \om_0 = \der x \wedge \der u + \der y \wedge \der v + \der z \wedge \der w.
\end{equation*}
In this setting, $T^* S^2$ is a closed symplectic submanifold of $T^* \R^3$ determined by the equations
\begin{align*}
  x^2 + y^2 + z^2 &= 1, \\
  xu + yv + zw &= 0,
\end{align*}
with the inherited symplectic form $\om$.

For convenience in our calculation, we also define a coordinate atlas on $T^* S^2$ by spelling out a set of local coordinates on several regions.

In $U_N = \Set{(x, y, z, u, v, w) \mmid z > 0}$, a neighborhood of stationary at the north pole $N = (0, 0, 1) \in S^2$, we use the polar coordinates of $(x, y)$ and $(u, v)$ to define local coordinates $\varphi_N = (\rho, \eta, \theta, \phi)$ determined by the parameterizations 
\begin{align*}
  (x, y, z) &= \Pa{\rho \cos \theta, \rho \sin \theta, \sqrt{1 - \rho^2}}, \\
  (u, v, w) &= \Pa{\eta \cos \phi, \eta \sin \phi, \frac{-\rho \eta \cos \delta}{\sqrt{1 - \rho^2}}},
\end{align*}
where $\delta = \phi - \theta$.
The range of these coordinates is $\varphi_N (U_N) = [0, 1) \times [0, \infty) \times (\R / 2\pi \Z)^2$; that is, we treat $\theta, \phi$ as $2\pi$-periodic coordinates, and ignore the non-injectivity of these coordinates at $(x, y) = 0$ or $(u, v) = 0$.
The symplectic form $\om_N$ on $U_N$ is the pullback of $\om_0$ by $\varphi_N$ as
\begin{multline*}
  \om_N = \varphi_N^* \om_0 = \frac{\cos \delta}{1 - \rho^2} \der \rho \wedge \der \eta + \frac{\rho^2 \eta \sin \delta}{1 - \rho^2} \der \rho \wedge \der \theta - \frac{\eta \sin \delta}{1 - \rho^2} \der \rho \wedge \der \phi \\
  - \rho \sin \delta \der \eta \wedge \der \theta + \rho \eta \cos \delta \der \theta \wedge \der \phi.
\end{multline*}

In $U_S = \Set{(x, y, z, u, v, w) \mmid z < 0}$, a neighborhood of stationary at the south pole $S = (0, 0, -1) \in S^2$, we use the polar coordinates of $(x, y)$ and $(u, v)$ to define local coordinates $\varphi_S = (\rho, \eta, \theta, \phi)$ determined by the parameterizations 
\begin{align*}
  (x, y, z) &= \Pa{\rho \cos \theta, \rho \sin \theta, -\sqrt{1 - \rho^2}}, \\
  (u, v, w) &= \Pa{\eta \cos \phi, \eta \sin \phi, \frac{\rho \eta \cos \delta}{\sqrt{1 - \rho^2}}}.
\end{align*}
The range of these coordinates is $\varphi_S (U_S) = [0, 1) \times [0, \infty) \times (\R / 2\pi \Z)^2$; that is, we treat $\theta, \phi$ as $2\pi$-periodic coordinates, and ignore the non-injectivity of these coordinates at $(x, y) = 0$ or $(u, v) = 0$.
The symplectic form $\om_S$ on $U_S$ is the pullback of $\om_0$ by $\varphi_S$ as
\begin{multline*}
  \om_S = \varphi_S^* \om_0 = \frac{\cos \delta}{1 - \rho^2} \der \rho \wedge \der \eta + \frac{\rho^2 \eta \sin \delta}{1 - \rho^2} \der \rho \wedge \der \theta - \frac{\eta \sin \delta}{1 - \rho^2} \der \rho \wedge \der \phi \\
  - \rho \sin \delta \der \eta \wedge \der \theta + \rho \eta \cos \delta \der \theta \wedge \der \phi.
\end{multline*}

In $U_E = \Set{(x, y, z, u, v, w) \mmid -1 < z < 1}$, a neighborhood of stationary on the equator $E = (\R^2 \times \Set{0}) \cap S^2$, we need to utilize the coordinates $w$ to define local coordinates $\varphi_E = (z, w, \theta, \phi)$ determined by the parameterizations
\begin{align*}
  (x, y, z) &= \Pa{\sqrt{1 - z^2} \cos \theta, \sqrt{1 - z^2} \sin \theta, z}, \\
  (u, v, w) &= \Pa{- \frac{z w}{\sqrt{1 - z^2}} \sec \delta \cos \phi, - \frac{z w}{\sqrt{1 - z^2}} \sec \delta \sin \phi, w}.
\end{align*}
The range of these coordinates is $\varphi_E (U_E) = (-1, 1) \times [0, \infty) \times (\R / 2\pi \Z)^2$; that is, we treat $\theta, \phi$ as $2\pi$-periodic coordinates, and ignore the non-injectivity of these coordinates at $(x, y) = 0$ or $(u, v) = 0$.
The symplectic form $\om_E$ on $U_E$ is the pullback of $\om_0$ by $\varphi_E$ as
\begin{multline*}
  \om_E = \varphi_E^* \om_0 = \frac{1}{1 - z^2} \der z \wedge \der w + \frac{1 + z^2}{1 - z^2} w \tan \delta \der z \wedge \der \theta 
  + z \tan \delta \der w \wedge \der \theta - z w \sec^2 \delta \der \theta \wedge \der \phi.
\end{multline*}

\section{Hamiltonian dynamics of a general spherical pendulum} \label{sec:spdlm-gen}

\subsection{The pendulum in local coordinates} \label{ssec:coord-spdlm-gen}

The \emph{spherical pendulum} with potential $V$ models the motion of a particle on the frictionless surface of the unit ball placed in a conservative force field with \emph{potential function} $V$.
In this system, where the external force has a potential, the \emph{mechanical energy} $H$ is always conserved.
If the force field is rotationally symmetric about the $z$-axis, or namely, the potential only depends on the $z$-coordinate, then the $z$-\emph{angular momentum} $J$ is also conserved.
In this case, we obtain a first integral $J$ for the Hamiltonian function $H$, and therefore, the \emph{momentum map} $F = (J, H)$ defines an \emph{integrable system} on $(T^* S^2, \om)$.

To write in coordinates, let $V \colon [-1, 1] \to \R$ be a smooth function.
We write out the functions $J, H \colon T^* S^2 \to \R$ as
\begin{equation} \begin{split} \label{eq:JH-in-cart-spdlm-gen}
  J &= xv - yu, \\
  H &= \frac12 (u^2 + v^2 + w^2) + V (z),
\end{split} \end{equation}
where one verifies that $\{J, H\} = 0$.

In the local coordinates $\varphi_N = (\rho, \eta, \theta, \phi)$, the momentum map reads
\begin{equation} \begin{split} \label{eq:JH-in-varphiN-spdlm-gen}
  J &= \rho \eta \sin \delta, \\
  H &= \frac{\eta^2}{2} \frac{1 - \rho^2 \sin^2 \delta}{1 - \rho^2} + V \Pa{\sqrt{1 - \rho^2}},
\end{split} \end{equation}
where $\delta = \phi - \theta$.

In the local coordinates $\varphi_S = (\rho, \eta, \theta, \phi)$, the momentum map reads
\begin{equation} \begin{split} \label{eq:JH-in-varphiS-spdlm-gen}
  J &= \rho \eta \sin \delta, \\
  H &= \frac{\eta^2}{2} \frac{1 - \rho^2 \sin^2 \delta}{1 - \rho^2} + V \Pa{-\sqrt{1 - \rho^2}}.
\end{split} \end{equation}

In the local coordinates $\varphi_E = (z, w, \theta, \phi)$, the momentum map reads
\begin{align*}
  J &= z w \tan \delta, \\
  H &= \frac{w^2}{2} \Pa{1 + \frac{z^2}{1 - z^2} \sec^2 \delta} + V (z).
\end{align*}

\subsection{The general behavior} \label{ssec:behavior-spdlm-gen}

From now on, we assume that $V \colon [-1, 1] \to [0, 1]$ is a unimodal function with maximum $V(-1) = V(1) = 1$ and minimum $V(0) = 0$, strictly decreasing in $[-1, 0]$ and strictly increasing in $[0, 1]$.

Since $J$ is related to a $2\pi$-periodic rotational symmetry \emph{\`a la} Noether, the system $(T^* S^2, \om, F)$ is semitoric except that $J$ is not proper.
The momentum map $F$ is proper, and so the fibers are compact.
There is one focus-focus singular value $(0, 1)$, whose fiber consists of two fixed points $(0, 0, \pm1, 0, 0, 0)$ and two orbits diffeomorphic to $S^1 \times \R$.
This fiber is a doubly pinched two-torus.

\begin{figure}[htb]
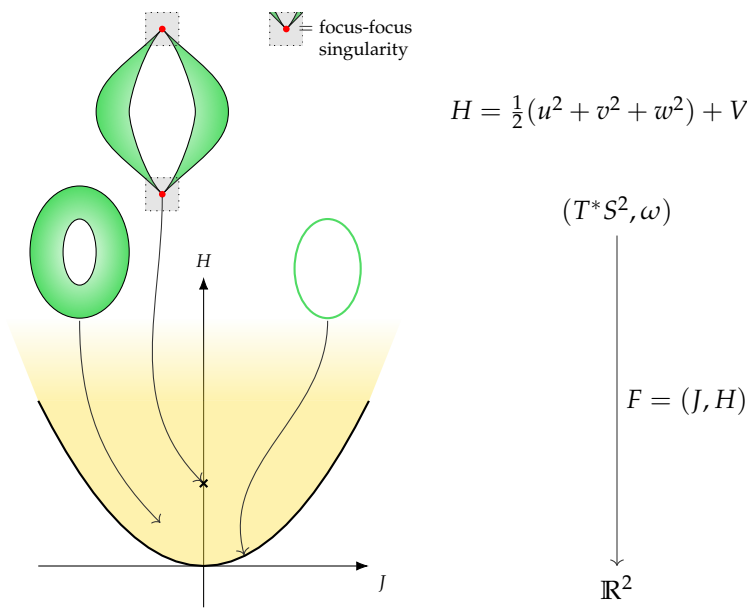

  \centering
  \inputfigure{momentum-image-2}
  \caption{The momentum map image and fibers.}
  \label{fig:momentum-image-2}
\end{figure}

To find the elliptic-regular singularities, one solves the minimum mechanical energy given a fixed angular momentum.
In fact, for $(j, h) = (J, H) (x, y, z, u, v, w)$ we have
\begin{align*}
  h &= \frac12 (u^2 + v^2 + w^2) + V (z) \geq \frac{j^2}{2(1 - z^2)} + V (z) \geq \frac{j^2}{2},
\end{align*}
by Cauchy--Schwarz inequality $(x^2 + y^2) (u^2 + v^2) \geq (xv - yu)^2$.
Then the momentum map image is the epigraph
\begin{align} \label{eq:JH-image-spdlm-gen}
  F (T^* S^2) = \Set{(j, h) \in \R^2 \mmid j \in \R, h \geq \frac{j^2}{2}},
\end{align}
and the fiber at each value $(j, \frac{j^2}{2})$ is a circle consisting of elliptic-regular points.
All the other fibers are regular two-tori; see \cref{fig:momentum-image-2}

There are action-angle coordinates in the preimage of a simply-connected region of regular values.
We are going to give a full calculation for the symmetric quadratic case in \cref{sec:action-angle}.

\subsection{The Hamiltonian dynamics} \label{ssec:ham-spdlm-gen}

Within the coordinate chart $(U_N, \varphi_N = (\rho, \eta, \theta, \phi))$, we differentiate the momentum map \cref{eq:JH-in-varphiN-spdlm-gen} to get 
\begin{equation*} \begin{split}
  \der J &= \eta \sin \delta \der \rho + \rho \sin \delta \der \eta - \rho \eta \cos \delta \der \theta + \rho \eta \cos \delta \der \phi, \\
  \der H &= \Pa{\frac{\rho \eta^2 (1 - \sin^2 \delta)}{(1 - \rho^2)^2} + \wt{V}_N' (\rho)} \der \rho + \eta \frac{1 - \rho^2 \sin^2 \delta}{1 - \rho^2} \der \eta + \frac{\rho^2 \eta^2 \sin \delta \cos \delta}{1 - \rho^2} (\der \theta - \der \phi),
\end{split} \end{equation*}
where $\delta = \phi - \theta$ and $\wt{V}_N (\rho) = V \Pa{\sqrt{1 - \rho^2}}$.
Then, the Hamiltonian vector fields associated with $J$ and $H$ are represented by
\begin{align*}
  X_J &= \om^{-1} \der J = \frac{\partial}{\partial \theta} + \frac{\partial}{\partial \phi}, \\
  X_H &= \om^{-1} \der H = \eta \cos \delta \Pa{\frac{\partial}{\partial \rho} - \Pa{\frac{\rho \eta (1 - \rho^2 \sin^2 \delta)}{1 - \rho^2} + \frac{(1 - \rho^2) \wt{V}_N' (\rho)}{\eta}} \,\frac{\partial}{\partial \eta}} \\
  &\quad + \eta \sin \delta \Pa{\frac{1}{\rho}\frac{\partial}{\partial \theta} + \Pa{\frac{\rho \eta (1 - \rho^2 \sin^2 \delta)}{1 - \rho^2} + \frac{(1 - \rho^2) \wt{V}_N' (\rho)}{\eta}} \frac{1}{\eta} \frac{\partial}{\partial \phi}}.
\end{align*}
The Hamiltonian equations with respect to $H$ are now ordinary differential equations as
\begin{align}
  \frac{\der \rho}{\der t} &= \eta \cos \delta, \label{eq:ham-eq-H-spdlm-gen-rho} \\
  \frac{\der \eta}{\der t} &= -\Pa{\frac{\rho \eta (1 - \rho^2 \sin^2 \delta)}{1 - \rho^2} + \frac{(1 - \rho^2) \wt{V}_N' (\rho)}{\eta}} \eta \cos \delta, \label{eq:ham-eq-H-spdlm-gen-eta} \\
  \frac{\der \theta}{\der t} &= \frac{\eta}{\rho} \sin \delta, \label{eq:ham-eq-H-spdlm-gen-theta} \\
  \frac{\der \phi}{\der t} &= \Pa{\frac{\rho \eta (1 - \rho^2 \sin^2 \delta)}{1 - \rho^2} + \frac{(1 - \rho^2) \wt{V}_N' (\rho)}{\eta}} \sin \delta. \label{eq:ham-eq-H-spdlm-gen-phi}
\end{align}

Within the coordinate chart $(U_S, \varphi_S = (\rho, \eta, \theta, \phi))$, we differentiate the momentum map \cref{eq:JH-in-varphiS-spdlm-gen} to get 
\begin{equation*} \begin{split}
  \der J &= \eta \sin \delta \der \rho + \rho \sin \delta \der \eta - \rho \eta \cos \delta \der \theta + \rho \eta \cos \delta \der \phi, \\
  \der H &= \Pa{\frac{\rho \eta^2 (1 - \sin^2 \delta)}{(1 - \rho^2)^2} + \wt{V}_S (\rho)} \der \rho + \eta \frac{1 - \rho^2 \sin^2 \delta}{1 - \rho^2} \der \eta + \frac{\rho^2 \eta^2 \sin \delta \cos \delta}{1 - \rho^2} (\der \theta - \der \phi),
\end{split} \end{equation*}
where $\delta = \phi - \theta$ and $\wt{V}_S (\rho) = V \Pa{-\sqrt{1 - \rho^2}}$.
Then, the Hamiltonian vector fields associated with $J$ and $H$ are represented by
\begin{align*}
  X_J &= \om^{-1} \der J = \frac{\partial}{\partial \theta} + \frac{\partial}{\partial \phi}, \\
  X_H &= \om^{-1} \der H = \eta \cos \delta \Pa{\frac{\partial}{\partial \rho} - \Pa{\frac{\rho \eta (1 - \rho^2 \sin^2 \delta)}{1 - \rho^2} + \frac{(1 - \rho^2) \wt{V}_S (\rho)}{\eta}} \,\frac{\partial}{\partial \eta}} \\
  &\quad + \eta \sin \delta \Pa{\frac{1}{\rho}\frac{\partial}{\partial \theta} + \Pa{\frac{\rho \eta (1 - \rho^2 \sin^2 \delta)}{1 - \rho^2} + \frac{(1 - \rho^2) \wt{V}_S (\rho)}{\eta}} \frac{1}{\eta} \frac{\partial}{\partial \phi}}.
\end{align*}
The Hamiltonian equations with respect to $H$ are now ordinary differential equations as
\begin{align*}
  \frac{\der \rho}{\der t} &= \eta \cos \delta, \\
  \frac{\der \eta}{\der t} &= -\Pa{\frac{\rho \eta (1 - \rho^2 \sin^2 \delta)}{1 - \rho^2} + \frac{(1 - \rho^2) \wt{V}_S (\rho)}{\eta}} \eta \cos \delta, \\
  \frac{\der \theta}{\der t} &= \frac{\eta}{\rho} \sin \delta, \\
  \frac{\der \phi}{\der t} &= \Pa{\frac{\rho \eta (1 - \rho^2 \sin^2 \delta)}{1 - \rho^2} + \frac{(1 - \rho^2) \wt{V}_S (\rho)}{\eta}} \sin \delta.
\end{align*}

\subsection{The explicit solutions} \label{ssec:sol-spdlm-gen}

One advantage of integrable systems is that their Hamiltonian equations can be explicitly solved by quadratures.
Since the trajectories are contained in the level sets of the momentum maps, we study their geometries.
Let $(j, h)$ be an arbitrary value in the momentum map image $F (T^* S^2)$.
The fiber is given by equations $J = j, H = h$.
Within the coordinate chart $(U_N, \varphi_N = (\rho, \eta, \theta, \phi))$, we express $\eta$ on this fiber in terms of $\rho$ via
\begin{equation} \label{eq:eta-rho-fiber-spdlm-gen}
  \eta^2 = j^2 + 2(1 - \rho^2) (h - \wt{V}_N (\rho)).
\end{equation}

Since $\rho \eta \sin \delta = j$ is constant on a trajectory, we see that $\sin \delta$ does not change sign.
By the continuity of $\delta$ in $t$ we now work on a shorter trajectory where $\epsilon = \sgn \cos \delta$ is also unchanged.
Substituting $\rho \eta \sin \delta = j$ into \cref{eq:ham-eq-H-spdlm-gen-rho}, and representing $\eta$ as a function of $\rho$ by \cref{eq:eta-rho-fiber-spdlm-gen}, we obtain that
\begin{equation*}
  \frac{\der \rho}{\der t} = \epsilon \frac{\sqrt{(1 - \rho^2) (2\rho^2 (h - \wt{V}_N (\rho)) - j^2)}}{\rho}.
\end{equation*}
This is an autonomous first-order differential equation separable in $\rho$ and can be solved for $\rho$ using integration.
If the initial value is given by $\rho_0 = \rho (0)$, then $\rho (t) = R_\epsilon (t; \rho_0)$, where $R_\epsilon$ are functions depending on $V$ and the values of $j, h$ presented implicitly as
\begin{align*}
  \int_{\rho_0}^{R_\epsilon (t; \rho_0)} \frac{\rho}{\sqrt{(1 - \rho^2)(2\rho^2(h - \wt{V}_N (\rho)) - j^2)}} \der \rho = \epsilon t.
\end{align*}

We find $\eta(t)$ immediately by \cref{eq:eta-rho-fiber-spdlm-gen} as
\begin{equation*}
  \eta (t) = \sqrt{j^2 + 2 (1 - R_\epsilon (t; \rho_0)^2) (h - \wt{V}_N \circ R_\epsilon (t; \rho_0))}.
\end{equation*}

Substituting $\rho (t) = R_\epsilon (t; \rho_0)$ into \cref{eq:ham-eq-H-spdlm-gen-rho}, we obtain that
\begin{align*}
  \frac{\der \theta}{\der \rho} &= \frac{\der \theta}{\der t} \Big/ \frac{\der \rho}{\der t} = \frac{\epsilon j}{\rho\sqrt{(1 - \rho^2) (2\rho^2 (h - \wt{V}_N (\rho)) - j^2)}}.
\end{align*}
This is an autonomous first-order differential equation separable in $\rho$ and can be solved for $\theta$ using integration.
If the initial value is given by $\theta_0 = \theta(0)$, and then
\begin{align*}
  \theta (t) = \theta_0 + \epsilon j \int_{\rho_0}^{R_\epsilon (t; \rho_0)} \frac{\der \rho}{\rho\sqrt{(1 - \rho^2) (2\rho^2 (h - \wt{V}_N (\rho)) - j^2)}}.
\end{align*}

Now, using the fact that
\begin{align*}
  \sin \delta (t) = \frac{j}{\rho (t) \eta (t)} = \frac{j}{R_\epsilon (t; \rho_0) \sqrt{j^2 + 2 (1 - R_\epsilon (t; \rho_0)^2) (h - \wt{V}_N \circ R_\epsilon (t; \rho_0))}},
\end{align*}
we obtain that, with the initial value $\phi_0 = \phi(0)$,
\begin{align*}
  \phi (t) &= \theta_0 + \arcsin \frac{j}{R_\epsilon (t; \rho_0) \sqrt{j^2 + 2 (1 - R_\epsilon (t; \rho_0)^2) (h - \wt{V}_N \circ R_\epsilon (t; \rho_0))}} \\
  &\quad + \epsilon j \int_{\rho_0}^{R_\epsilon (t; \rho_0)} \frac{\der \rho}{\rho\sqrt{(1 - \rho^2) (2\rho^2 (h - \wt{V}_N (\rho)) - j^2)}}.
\end{align*}

Denote by $\Phi_H^t \in \operatorname{Ham} (T^* S^2, \om)$, $t \in \R$ the flow of $X_H$.
Summarizing the above formulas, we solve the Hamiltonian equations for $H$ as
\begin{equation*} \begin{split}
  \rho (t) &= R_\epsilon (t; \rho_0), \\
  \eta (t) &= \sqrt{j^2 + 2 (1 - R_\epsilon (t; \rho_0)^2) (h - \wt{V}_N \circ R_\epsilon (t; \rho_0))}, \\
  \theta (t) &= \theta_0 + \epsilon j \int_{\rho_0}^{R_\epsilon (t; \rho_0)} \frac{\der \rho}{\rho\sqrt{(1 - \rho^2) (2\rho^2 (h - \wt{V}_N (\rho)) - j^2)}}, \\
  \phi (t) &= \phi_0 - \arcsin \frac{j}{\rho_0 \eta_0} + \arcsin \frac{j}{R_\epsilon (t; \rho_0) \sqrt{j^2 + 2 (1 - R_\epsilon (t; \rho_0)^2) (h - \wt{V}_N \circ R_\epsilon (t; \rho_0))}} \\
  &\quad + \epsilon j \int_{\rho_0}^{R_\epsilon (t; \rho_0)} \frac{\der \rho}{\rho\sqrt{(1 - \rho^2) (2\rho^2 (h - \wt{V}_N (\rho)) - j^2)}},
\end{split} \end{equation*}
with initial values $(\rho, \eta, \theta, \phi) (0) = (\rho_0, \eta_0, \theta_0, \phi_0)$ and $\epsilon = \sgn \cos \delta$ being constant.
The Hamiltonian flow for $H$ is given by $(\rho, \eta, \theta, \phi) (t) = \Phi_H^t (\rho_0, \eta_0, \theta_0, \phi_0)$.
If $\cos \delta$ does change sign, it does so only at some $t = t_1$ when $\delta (t_1) = \pm \frac\pi2$.
Then, we can determine the sign of 
\begin{align*}
  \frac{\der \delta}{\der t} (t_1) &= \frac{(1 - \rho^2) (\rho \wt{V}_N' (\rho) - \eta^2)}{\rho \eta}.
\end{align*}

Denote by $\Phi_J^s \in \operatorname{Ham} (T^* S^2, \om)$, $s \in \R$ the flow of $X_J$.
The Hamiltonian equations for $J$ are fairly simple as
\begin{align*}
  \frac{\der \rho}{\der s} &= 0, &\frac{\der \eta}{\der s} &= 0, &\frac{\der \theta}{\der s} &= 1, &\frac{\der \phi}{\der s} &= 1,
\end{align*}
whose solutions are given by
\begin{align*}
  \rho (s) &= \rho_0, &\eta (s) &= \eta_0, &\theta (s) &= \theta_0 + s, &\phi (s) &= \phi_0 + s,
\end{align*}
with the initial condition $(\rho, \eta, \theta, \phi) (0) = (\rho_0, \eta_0, \theta_0, \phi_0)$.
The Hamiltonian flow for $J$ is given by $(\rho, \eta, \theta, \phi) (s) = \Phi_J^s (\rho_0, \eta_0, \theta_0, \phi_0)$.

To evaluate the joint flow $(\rho, \eta, \theta, \phi) (s, t) = \Phi_J^s \circ \Phi_H^t (\rho_0, \eta_0, \theta_0, \phi_0)$, we observe that $\delta$ is constant under the flow of $X_J$.
Composing the two single flows, we have that
\begin{equation} \begin{split} \label{eq:sol-spdlm-gen}
  \rho (s, t) &= R_\epsilon (t; \rho_0), \\
  \eta (s, t) &= \sqrt{j^2 + 2 (1 - R_\epsilon (t; \rho_0)^2) (h - \wt{V}_N \circ R_\epsilon (t; \rho_0))}, \\
  \theta (s, t) &= \theta_0 + \epsilon j \int_{\rho_0}^{R_\epsilon (t; \rho_0)} \frac{\der \rho}{\rho\sqrt{(1 - \rho^2) (2\rho^2 (h - \wt{V}_N (\rho)) - j^2)}} + s, \\
  \phi (s, t) &= \phi_0 - \arcsin \frac{j}{\rho_0 \eta_0} + \arcsin \frac{j}{R_\epsilon (t; \rho_0) \sqrt{j^2 + 2 (1 - R_\epsilon (t; \rho_0)^2) (h - \wt{V}_N \circ R_\epsilon (t; \rho_0))}} \\
  &\quad + \epsilon j \int_{\rho_0}^{R_\epsilon (t; \rho_0)} \frac{\der \rho}{\rho\sqrt{(1 - \rho^2) (2\rho^2 (h - \wt{V}_N (\rho)) - j^2)}} + s.
\end{split} \end{equation}

We know that the regular fiber given by $F = (j, h)$ is diffeomorphic to $\mathbb{T}^2$ and the flows of $X_J$ and $X_H$ are quasi-periodic.
With the explicit formulas \cref{eq:sol-spdlm-gen} for the joint flow, one could calculate the period lattice so as to find the action-angle coordinates.
However, the calculation could not be fully done in a single chart, and one should write down the solutions carefully under the other coordinates and keep track of the coordinate transformations.
For simplicity, we only give a full answer to those with symmetric quadratic potentials.

\section{The spherical pendulum with a symmetric quadratic potential and its dynamics} \label{sec:spdlm-quad}

In this section, we focus on the spherical pendulum with a \emph{symmetric quadratic potential} function $V (z) = z^2$ and solve the Hamiltonian equations explicitly.

\subsection{The Hamiltonian dynamics} \label{ssec:ham-spdlm-quad}

Specialized in the case $V (z) = z^2$, the momentum map $F = (J, H)$ is written
\begin{equation} \begin{split} \label{eq:JH-in-cart-spdlm-quad}
  J &= xv - yu, \\
  H &= \frac12 (u^2 + v^2 + w^2) + z^2. 
\end{split} \end{equation}

In the local coordinates $\varphi_N = (\rho, \eta, \theta, \phi)$, the momentum map reads
\begin{align*}
  J &= \rho \eta \sin \delta, \\
  H &= \frac{\eta^2}{2} \frac{1 - \rho^2 \sin^2 \delta}{1 - \rho^2} + 1 - \rho^2,
\end{align*}
where $\delta = \phi - \theta$.

Computing the differentials, we obtain
\begin{align*}
  \der J &= \eta \sin \delta \der \rho + \rho \sin \delta \der \eta - \rho \eta \cos \delta \der \theta + \rho \eta \cos \delta \der \phi, \\
  \der H &= \Pa{\rho \eta^2 \frac{\cos^2 \delta}{(1 - \rho^2)^2} - 2 \rho} \der \rho + \eta \frac{1 - \rho^2 \sin^2 \delta}{1 - \rho^2} \der \eta + \frac{\rho^2 \eta^2 \sin \delta \cos \delta}{1 - \rho^2} (\der \theta - \der \phi).
\end{align*}
Then, the Hamiltonian vector fields associated with $J$ and $H$ are represented by
\begin{align*}
  X_J &= \om^{-1} \der J = \frac{\partial}{\partial \theta} + \frac{\partial}{\partial \phi}, \\
  X_H &= \om^{-1} \der H = \eta \cos \delta \Pa{\frac{\partial}{\partial \rho} - \Pa{\frac{\rho \eta (1 - \rho^2 \sin^2 \delta)}{1 - \rho^2} - \frac{2\rho (1 - \rho^2)}{\eta}} \,\frac{\partial}{\partial \eta}} \\
  &\quad + \eta \sin \delta \Pa{\frac{1}{\rho}\frac{\partial}{\partial \theta} + \Pa{\frac{\rho \eta (1 - \rho^2 \sin^2 \delta)}{1 - \rho^2} - \frac{2\rho (1 - \rho^2)}{\eta}} \frac{1}{\eta} \frac{\partial}{\partial \phi}}.
\end{align*}
The ordinary differential equations for the trajectory of $X_H$ are now
\begin{align*}
  \frac{\der \rho}{\der t} &= \eta \cos \delta, \\
  \frac{\der \eta}{\der t} &= -\rho \eta \cos \delta \Pa{\frac{\eta (1 - \rho^2 \sin^2 \delta)}{1 - \rho^2} - \frac{2 (1 - \rho^2)}{\eta}}, \\
  \frac{\der \theta}{\der t} &= \frac{\eta}{\rho} \sin \delta, \\
  \frac{\der \phi}{\der t} &= \rho \sin \delta \Pa{\frac{\eta (1 - \rho^2 \sin^2 \delta)}{1 - \rho^2} - \frac{2 (1 - \rho^2)}{\eta}}.
\end{align*}

\subsection{The explicit solutions} \label{ssec:sol-spdlm-quad}

In this case, we are able to write the solutions \cref{eq:sol-spdlm-gen} in a closed form in terms of the elliptic integrals.
For the sake of brevity, we write
\begin{equation}  \label{eq:def-k-n}
  k = \sqrt{\frac{1 + h - \sqrt{(1 - h)^2 + 2j^2}}{1 + h + \sqrt{(1 - h)^2 + 2j^2}}}, \qquad
  n = \frac{1 + h - \sqrt{(1 - h)^2 + 2j^2}}{2}.
\end{equation}

First, we integrate
\begin{align*}
  &\quad \int \frac{\rho}{\sqrt{(1 - \rho^2) \Pa{2\rho^2 (h + \rho^2 - 1) - j^2}}} \der \rho \\
  &= -\frac{k}{\sqrt{2 n}} \int \frac{1}{\sqrt{1 - k^2 \sin^2 \gamma}} \der \gamma = -\frac{k}{\sqrt{2 n}} \,F \Pa{\cdot, k} + \text{const},
\end{align*}
where we apply the change of variables $\gamma = \arcsin \frac{\sqrt{1 - \rho^2}}{\sqrt{n}}$ and $F$ is the elliptic integral of the first kind as in \cref{eq:ell-int-1}.

If the initial value is given by $\rho_0 = \rho (0)$, and then $\rho = \rho (t)$ are presented implicitly as
\begin{align*}
  F \Pa{\arcsin \frac{\sqrt{1 - \rho^2}}{\sqrt{n}}, k} - F \Pa{\arcsin \frac{\sqrt{1 - \rho_0^2}}{\sqrt{n}}, k} = -\frac{\epsilon \sqrt{2 n}}{k} t.
\end{align*}
The function $R_\epsilon$ for $V (z) = z^2$ is presented explicitly as
\begin{align*}  
  R_\epsilon (t; \rho_0) = \sqrt{1 - n \operatorname{sn}^2 \Pa{F \Pa{\arcsin \frac{\sqrt{1 - \rho_0^2}}{\sqrt{n}}, k} - \frac{\epsilon \sqrt{2 n}}{k} t, k}},
\end{align*}
where $\operatorname{sn}$ is the Jacobi elliptic sine function as in \cref{eq:jacobi-sn}.

Second, we integrate
\begin{align*}
  &\quad \int \frac{\der \rho}{\rho\sqrt{(1 - \rho^2) (2\rho^2 (h + \rho^2 - 1) - j^2)}} \\
  &= -\frac{k}{\sqrt{2 n}} \int \frac {1}{(1 - n \sin^2 \gamma) \sqrt{1 - k^2 \sin^2 \gamma}} \der \gamma = -\frac{k}{\sqrt{2 n}} \,\Pi \Pa{\cdot, n, k} + \text{const},
\end{align*}
where we apply the change of variables $\gamma = \arcsin \frac{\sqrt{1 - \rho^2}}{\sqrt{n}}$ and $\Pi$ is the elliptic integral of the third kind as in \cref{eq:ell-int-3}.

If the initial value is given by $\rho_0 = \rho (0)$, and then $\theta = \theta (t)$ are presented explicitly as
\begin{align*}
  &\int_{\rho_0}^{R_\epsilon (t; \rho_0)} \frac{\der \rho}{\rho\sqrt{(1 - \rho^2) (2\rho^2 (h + \rho^2 - 1) - j^2)}} \\
  &= -\frac{k}{\sqrt{2 n}} \,\Pi \Pa{\arcsin \frac{\sqrt{1 - \rho^2}}{\sqrt{n}}, n, k} + \frac{k}{\sqrt{2 n}} \,\Pi \Pa{\arcsin \frac{\sqrt{1 - \rho_0^2}}{\sqrt{n}}, n, k}.
\end{align*}

To simplify our formulas, let
\begin{align*}
  \gamma (t) &= \operatorname{am} \Pa{F \Pa{\gamma_0, k} - \frac{\epsilon \sqrt{2 n}}{k} t, k}, & \gamma_0 &= \arcsin \frac{\sqrt{1 - \rho_0^2}}{\sqrt{n}},
\end{align*}
where $\operatorname{am}$ is the Jacobi amplitude function as in \cref{eq:jacobi-am}.
The joint flow $(\rho, \eta, \theta, \phi) (s, t) = \Phi_J^s \circ \Phi_H^t (\rho_0, \eta_0, \theta_0, \phi_0)$ has an explicit form as
\begin{equation} \begin{split} \label{eq:sol-spdlm-quad}
  \rho (s, t) &= \sqrt{1 - n \sin^2 \gamma(t)}, \\
  \eta (s, t) &= \sqrt{j^2 + 2 h n \sin^2 \gamma(t) - 2 n^2 \sin^4 \gamma(t)}, \\
  \theta (s, t) &= \theta_0 + \epsilon \sqrt{2 n}\, j k \,\Pi \Pa{\gamma_0, n, k} - \epsilon \sqrt{2n}\, j k \,\Pi \Pa{\gamma (t), n, k} + s, \\
  \phi (s, t) &= \phi_0 + \arcsin \frac{j}{\rho_0 \eta_0} + \epsilon \sqrt{2 n}\, j k \,\Pi \Pa{\gamma_0, n, k} - \epsilon \sqrt{2n}\, j k \,\Pi \Pa{\gamma (t), n, k} \\
  &\quad- \arcsin \frac{j}{\sqrt{j^2 + 2 n^3 k^{-2} \sin^2 \gamma(t) - 2 (1 + h) n^2 \sin^4 \gamma(t) + 2 n^3 \sin^6 \gamma(t)}} + s.
\end{split} \end{equation}

\section{Action-angle coordinates} \label{sec:action-angle}

\subsection{Introduction of action-angle coordinates} \label{ssec:intro-action-angle}

\emph{Action-angle coordinates} are a set of canonical coordinates used to describe integrable Hamiltonian systems.
They are particularly useful for studying periodic and quasi-periodic orbits. 
From Arnold's book \cite{MR997295}, we know that for an integrable system $(M, \om, F)$ we can find local coordinates
\begin{align*}
  (A_1, A_2, \alpha_1, \alpha_2) \colon F^{-1} (U) \to \R^2 \times (\R / 2\pi\Z)^2,
\end{align*}
where $U$ is a simply-connected open set of regular values of a proper momentum map $F$ such that
\begin{itemize}
  \item the momentum map $F$ depends only on the \emph{action coordinates} $A_1$ and $A_2$,
  \item the coordinates are symplectic, in the sense that $\omega = \der A_1 \wedge \der \alpha_1 + \der A_2 \wedge \der \alpha_2$, and
  \item the evolutions of the \emph{angle coordinates} $\alpha_1$ and $\alpha_2$ under Hamiltonian equations for $J$ and $H$ are linear in time.
\end{itemize}
The Hamiltonian vector fields for $A_1$ and $A_2$ are extremely simple:
\begin{align*}
  X_{A_1} = \frac{\partial}{\partial \alpha_1}, \qquad X_{A_2} = \frac{\partial}{\partial \alpha_2}.
\end{align*}

The action-angle coordinates are closely related to the periodic lattice: the Hamiltonian flows for $A_1$ and $A_2$ are both $2\pi$-periodic.
In the semitoric case, we always take the first action variable $A_1 = J$, since $J$ generates a Hamiltonian circle action.
As for the second action variable $A_2$, we define
\begin{equation} \label{eq:X-A2-in-JH}
  2\pi X_{A_2} = S (J, H) X_J + T (J, H) X_H,
\end{equation}
where $T$ denotes the time it takes a particle to return to the initial $X_J$-orbit under the Hamiltonian flow associated with $H$, and $S$ denotes the time it takes the particle to return to the initial status under the rotation by the Hamiltonian flow associated with $J$.
Both $S$ and $T$ depend only on the value of $F$, not the position on the fiber.

For the angle coordinates $\alpha_1$ and $\alpha_2$, since $\Set{\alpha_1 = \alpha_2 = 0}$ is Lagrangian, we need to specify a Lagrangian submanifold $L$ on which the angle coordinates vanish.
Then, we use the linearity of their evolutions under the Hamiltonian flows associated with $J$ and $H$ to extend the angle coordinates along the fiber.
The key mechanism is $\frac{\der \alpha_2}{\der t} = 1$ under the flow of $X_{A_2}$ as in \cref{eq:X-A2-in-JH}.

\subsection{Finding the periods} \label{ssec:period-spdlm-quad}

As in \cref{eq:JH-image-spdlm-gen} we have determined the momentum map image
\begin{align*}
  F (T^* S^2) = \Set{(j, h) \in \R^2 \mmid j \in \R, h \geq \frac{j^2}{2}},
\end{align*}
which consists of three strata: the curve $\partial F (T^* S^2) = \Set{j \in \R \mmid (j, \frac{j^2}{2})}$ of elliptic-regular values, the unique focus-focus value $(0, 1)$, and the region $B_\mathrm{r}$ of regular values.

For any value $(j, h)$ in
\begin{align*}
  B_\mathrm{r} = \Set{(j, h) \in \R^2 \mmid j \in \R, h > \frac{j^2}{2}, (j, h) \neq (0, 1)},
\end{align*}
there is a \emph{period lattice} $\Lambda_{j, h}$ in $T^*_{j, h} B_\mathrm{r}$.
Identifying the cotangent space with $\R^2$, we could write
\begin{align*}
  \Lambda_{j, h} = \Set{\frac{(s, t)}{2\pi} \in \R^2 \mmid \Res{\Phi_J^s \circ \Phi_H^t}_{F^{-1} (j, h)} = \identity}.
\end{align*}

In this subsection, we calculate the period lattice $\Lambda_{j, h}$ explicitly of $F = (J, H)$ for $V (z) = z^2$.
In fact, the calculation could not be fully done in one local chart $(U_N, \varphi_N)$, in which $\rho \in [0, 1)$ and $z \in (0, 1]$.
We should also be careful when we use the local coordinate representations as $\rho \to 1-$ and $z \to 0+$, by either transforming among coordinates or tracking the limit.
We fix a value $(j, h)$, and start with the point $(1, 0, 0)$ at the intersection of the positive $x$-axis and the equator, with the momentum $(0, j, \sqrt{2h - j^2})$.
Under these specific initial value conditions, we have $\rho_0 = 1$, $\eta_0 = j$, $\theta_0 = \gamma_0 = 0$, $\phi_0 = \delta_0 = \frac{\pi}{2}$, $\epsilon = -1$, and $\gamma (t) = \operatorname{am} \Pa{\frac{\sqrt{2 n}}{k} t, k}$.
The coordinates have rather neat expressions
\begin{equation} \begin{split} \label{eq:sol-rho-theta-L-spdlm-quad}
  \rho (s, t) &=\sqrt{1 - n \operatorname{sn}^2 \Pa{\frac{\sqrt{2 n}}{k} t, k}}, \\
  \theta (s, t) &= \frac{j k}{\sqrt{2 n}} \,\Pi \Pa{\operatorname{am} \Pa{\frac{\sqrt{2 n}}{k} t, k}, n, k} + s.
\end{split} \end{equation}
We would terminate this trajectory at $t = t_1$ when $\rho$ attains its minimum and $\cos \delta = 0$ at the first time.
Then, we rotate by $s = s_1$ so that $\theta (s_1, t_1) = \theta_0 = 0$.
At this moment, $\sin \delta (t_1) = 1$ since it does not change sign, and then by $\rho (t_1) \eta (t_1) = j$ and $\frac{\eta^2 (t_1)}{2} + 1 - \rho^2 (t_1) = h$ we find the only solution in $(0, 1)$ as
\begin{align*}
  \rho (s_1, t_1) = \sqrt{\frac{1 - h + \sqrt{(1 - h)^2 + 2j^2}}{2}}
\end{align*}
and then $\operatorname{am} \Pa{\frac{\sqrt{2 n}}{k} t_1, k} = \frac{\pi}{2}$, which implies that
\begin{align*}
  t_1 = \frac{k}{\sqrt{2 n}} \,F \Pa{\frac{\pi}{2}, k} = \frac{k}{\sqrt{2 n}} \,K \Pa{k},
\end{align*}
where $K$ is the complete elliptic integral of the first kind as in \cref{eq:ell-int-1-comp}.
At this moment, if $\theta (s_1, t_1) = \theta_0 = 0$ we have
\begin{align*}
  s_1 = -\frac{j k}{\sqrt{2 n}} \,\Pi \Pa{\frac{\pi}{2}, n, k} = -\frac{j k}{\sqrt{2 n}} \,\Pi \Pa{n, k},
\end{align*}
where $\Pi$ is the complete elliptic integral of the third kind as in \cref{eq:ell-int-3-comp}.

We shall benefit from the symmetry of the system.
If the trajectory continues, it takes another $t_1$ to return to the equator with a momentum pointing downwards.
After moving in the southern hemisphere, it returns to the equator again at $T = 4 t_1$ with the same value of $w$ as the beginning, and by rotating by $S = 4 s_1$ we arrive at the the starting point in $T^* S^2$.
Therefore, we obtain a generator of the periods as
\begin{align} \label{eq:def-S-T-spdlm-quad}
  (S, T) := (4 s_1, 4 t_1) = \Pa{-\frac{2 \sqrt2 j k}{\sqrt{n}} \,\Pi \Pa{n, k}, \frac{2 \sqrt2 k}{\sqrt{n}} \,K \Pa{k}}.
\end{align}

Note that we only calculated the periods \cref{eq:def-S-T-spdlm-quad} for $j \neq 0$ and $2h - j^2 > 0$.
For regular values with $j = 0$, we calculate the periods separately and realize that they fall into the same general formula.
So, \cref{eq:def-S-T-spdlm-quad} applies for all $(j, h) \in B_\mathrm{r}$.
We conclude that
\begin{align*}
  \Lambda_{j, h} = (1, 0) \Z + \Pa{-\frac{2 \sqrt2 k}{\sqrt{n}} \,\Pi \Pa{n, k}, \frac{2 \sqrt2 k}{\sqrt{n}} \,K \Pa{k}} \Z
\end{align*}
for $(j, h) \in B_\mathrm{r}$.
Do not forget an obvious generator as $(1, 0)$.

For convenience, we present the special cases below.
If $j = 0$ and $h < 1$, then
\begin{align*}
  \Lambda_{j, h} = (1, 0) \Z + \Pa{0, 2 \sqrt2 \,K \Pa{\sqrt{h}}} \Z.
\end{align*}
If $j = 0$ and $h > 1$, then
\begin{align*}
  \Lambda_{j, h} = (1, 0) \Z + \Pa{0, \frac{2 \sqrt2}{\sqrt{h}} \,K \Pa{\frac{1}{\sqrt{h}}}} \Z.
\end{align*}

For those boundary values with $2h - j^2 = 0$, the degenerate period lattice has rank one:
\begin{align*}
  \Lambda_{j, h} = (1, 0) \Z.
\end{align*}

\subsection{Finding the action-angle coordinates} \label{ssec:action-angle-spdlm-quad}

We have been prepared to compute the action coordinates $(A_1, A_2) \colon F^{-1} (B_{\mathrm{r}}) \to \R^2$.
Since $A_1 = J$, we are left to find $A_2$, which can be determined by \cref{eq:X-A2-in-JH}, as
\begin{align*}
  \der A_2 &= -\frac{\sqrt 2 J k}{\pi \sqrt{n}} \,\Pi \Pa{n, k} \der J + \frac{\sqrt 2 k}{\pi \sqrt{n}} \,K (k) \der H,
\end{align*}
for $2 H > J^2$ and either $J \neq 0$ or $0 < H < 1$.
Utilizing the change of variables in \cref{eq:def-k-n}, we write $\der A_2$ in terms of $k$ and $n$ for $0 < k \leq \sqrt{n} < 1$ as following
\begin{equation*}
  \der A_2 = \frac{\sqrt2}{\pi \sqrt{n}} \Pa{\frac{1 + k^2}{k} \,K (k) + \frac{2n - 1 - k^2}{k} \,\Pi (n,k)} \der n - \frac{2 \sqrt{2 n}}{\pi k^2} \Pa{K (k) + (n - 1) \Pi (n, k)} \der k.
\end{equation*}
By integration, we obtain
\begin{equation} \label{eq:A2-in-nk}
  A_2 =  \frac{2 \sqrt2}{\pi} \Pa{\frac{n - 1}{\sqrt{n}} k K(k) + \frac{\sqrt{n}}{k} E(k) - \frac{(n - 1)(k^2 - n)}{\sqrt{n} k} \Pi(n, k)}.
\end{equation}

Let $L = \Set{(1, 0, 0, 0, v, w) \in T^* S^2 \mmid w > 0}$ be the upper half of the Lagrangian fiber $T^*_{(1, 0, 0)} S^2$, on which we set the angle coordinates $(\alpha_1, \alpha_2) = (0, 0)$.
Given $(j, h) \in B_{\mathrm{r}}$, there is a unique point $(1, 0, 0, 0, j, \sqrt{2h - j^2}) \in L$ with $F = (j, h)$.
Recall $S$ and $T$ in \cref{eq:def-S-T-spdlm-quad}.
Any other point $(x, y, z, u, v, w)$ on $F^{-1} (j, h)$ is uniquely represented by $(s, t) \in (\R / S \Z) \times (\R / T \Z)$ according to
\begin{align*}
  \Phi_J^s \circ \Phi_H^t \Pa{1, 0, 0, 0, j, \sqrt{2h - j^2}} = (x, y, z, u, v, w),
\end{align*}
which can be solved by \cref{eq:sol-rho-theta-L-spdlm-quad} as
\begin{equation} \begin{split} \label{eq:s-t-formula-L-upward-spdlm-quad}
  s &= \theta - \frac{j k}{\sqrt{2 n}} \,\Pi \Pa{\arcsin \frac{z}{\sqrt{n}}, n, k}, \\
  t &= \frac{k}{\sqrt{2 n}} \,F \Pa{\arcsin \frac{z}{\sqrt{n}}, k}.
\end{split} \end{equation}
These formulas work for points with $0 < z < 1$ and $w > 0$.
By a limiting and a reflection arguments, we observe that \cref{eq:s-t-formula-L-upward-spdlm-quad} works in the case of $-1 \leq z \leq 1$ and $w \geq 0$, and when $w < 0$, we have
\begin{equation*} \begin{split}
  s &= \theta - \frac12 S + \frac{j k}{\sqrt{2 n}} \,\Pi \Pa{\arcsin \frac{z}{\sqrt{n}}, n, k}, \\
  t &= \frac12 T - \frac{k}{\sqrt{2 n}} \,F \Pa{\arcsin \frac{z}{\sqrt{n}}, k}.
\end{split} \end{equation*}

The parameters $(s, t)$ on the momentum map fiber are related to $(\alpha_1, \alpha_2)$ by the linear relation
\begin{align*}
  s \der J + t \der H &= \alpha_1 \der A_1 + \alpha_2 \der A_2 = \alpha_1 \der J + \alpha_2 \Pa{\frac{S}{2\pi} \der J + \frac{T}{2\pi} \der H}.
\end{align*}
Therefore, we get
\begin{align*}
  \alpha_1 &= s - S \frac{t}{T} = \begin{cases}
    \theta - \frac{J k}{\sqrt{2 n}} \Pa{\Pi \Pa{\arcsin \frac{z}{\sqrt{n}}, n, k} - \frac{F \Pa{\arcsin \frac{z}{\sqrt{n}}, k}}{K \Pa{k}} \Pi \Pa{n, k}} & \text{if }w \geq 0, \\
    \theta + \frac{J k}{\sqrt{2 n}} \Pa{2 \Pi \Pa{n, k} + \Pi \Pa{\arcsin \frac{z}{\sqrt{n}}, n, k} + \frac{F \Pa{\arcsin \frac{z}{\sqrt{n}}, k}}{K \Pa{k}} \Pi \Pa{n, k}} & \text{if }w < 0;
  \end{cases} \\
  \alpha_2 &= 2\pi \frac{t}{T} = \begin{cases}
    \frac{\pi}{2} \frac{F \Pa{\arcsin \frac{z}{\sqrt{n}}, k}}{K \Pa{k}} & \text{if }w \geq 0, \\
    \pi - \frac{\pi}{2} \frac{F \Pa{\arcsin \frac{z}{\sqrt{n}}, k}}{K \Pa{k}} & \text{if }w < 0.
  \end{cases}
\end{align*}

\begin{figure}[htb]
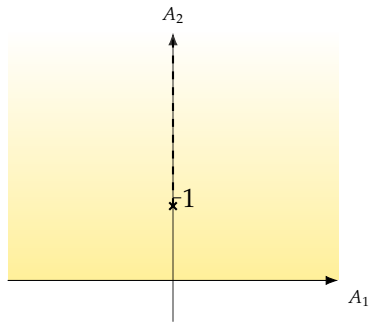

  \centering
  \inputfigure{action-coord}
  \caption{The image of the action coordinates $A$.}
  \label{fig:action-coord}
\end{figure}

Summarizing, the action-angle coordinates $(A_1, A_2, \alpha_1, \alpha_2)$ for the symmetric quadratic spherical pendulum for $0 < k \leq \sqrt{n} < 1$ are 
\begin{equation} \begin{split} \label{eq:action-angle-sq-spdlm}
  A_1 &= J, \\
  A_2 &=  \frac{2 \sqrt2}{\pi} \Pa{\frac{n - 1}{\sqrt{n}} k K(k) + \frac{\sqrt{n}}{k} E(k) - \frac{(n - 1)(k^2 - n)}{\sqrt{n} k} \Pi(n, k)}, \\
  \alpha_1 &= \begin{cases}
    \theta - \frac{J k}{\sqrt{2 n}} \Pa{\Pi \Pa{\arcsin \frac{z}{\sqrt{n}}, n, k} - \frac{F \Pa{\arcsin \frac{z}{\sqrt{n}}, k}}{K \Pa{k}} \Pi \Pa{n, k}} & \text{if }w \geq 0, \\
    \theta + \frac{J k}{\sqrt{2 n}} \Pa{2 \Pi \Pa{n, k} + \Pi \Pa{\arcsin \frac{z}{\sqrt{n}}, n, k} + \frac{F \Pa{\arcsin \frac{z}{\sqrt{n}}, k}}{K \Pa{k}} \Pi \Pa{n, k}} & \text{if }w < 0;
  \end{cases} \\
  \alpha_2 &= \begin{cases}
    \frac{\pi}{2} \frac{F \Pa{\arcsin \frac{z}{\sqrt{n}}, k}}{K \Pa{k}} & \text{if }w \geq 0, \\
    \pi - \frac{\pi}{2} \frac{F \Pa{\arcsin \frac{z}{\sqrt{n}}, k}}{K \Pa{k}} & \text{if }w < 0.
  \end{cases}
\end{split} \end{equation} 
where the angular momentum $J$ and the mechanical energy $H$ are defined in \cref{eq:JH-in-cart-spdlm-quad}, and the parameters are as in \cref{eq:def-k-n}.
Note that the coordinate $\theta$ is ambiguous at $(0, 0, \pm 1, u, v, 0)$; at these points we replace $\theta$ by $\phi$ in the formula for $\alpha_1$.
The image of $(A_1, A_2)$ and the focus-focus value are shown in \cref{fig:action-coord}.

\subsection{The monodromy of action coordinates}

As $J \to 0$ and for a fixed $H > 1$, we find
\begin{equation*}
  A_2 \to \frac{2 \sqrt2}{\pi k} E (k) = \frac{2 \sqrt{2 H}}{\pi} E \Pa{\frac{1}{\sqrt{H}}}.
\end{equation*}
Overall, $A_2$ can be extended to a continuous function on $T^* S^2$, which we still denote by $A_2$.
However, $A_2$ is not differentiable along $J = 0, H > 1$, caused by its monodromy.

Recall $\der A_2$ as in \cref{eq:A2-in-nk}.
As $J \to 0-$ and for $H > 1$, we calculate the limit
\begin{align*}
  \der A_2 &\to \der J + \frac{\sqrt{2}}{\pi \sqrt{H}} \,K \Pa{\frac{1}{\sqrt{H}}} \der H =: \der A_2^-.
\end{align*}
As $J \to 0+$ and for $H > 1$, we calculate the limit
\begin{align*}
  \der A_2 &\to -\der J + \frac{\sqrt{2}}{\pi \sqrt{H}} \,K \Pa{\frac{1}{\sqrt{H}}} \der H =: \der A_2^+.
\end{align*}
Consider a counter-clockwise loop based at $(0, h) \in F (T^* S^2)$ for some $h > 1$ which winds around $(0, 1)$ for one cycle.
Then $\der A$ changes from $\der A^- := (\der J, \der A_2^-)$ to $\der A^+ := (\der J, \der A_2^+)$ along this loop.
The \emph{monodromy} of the action coordinates along this loop is defined as
\begin{equation} \label{eq:monodromy-sq-spdlm}
  \der A^+ \circ (\der A^-)^{-1} = \begin{pmatrix}
    1 & 0 \\ 2 & 1                                               
  \end{pmatrix}.
\end{equation}
This coincides with the monodromy around focus-focus values; see, for instance, \cite{MR1937606}.
Moreover, according to \cref{eq:ell-int-2-comp-asym}, $\der A_2$ blows up to infinity in a logarithmic manner, as $(J, H) \to (0, 1)$.

As a comparison, when $J = 0$ and $0 < H < 1$, we find
\begin{align*}
  A_2 &= \frac{2 \sqrt2}{\pi} \Pa{(H - 1) K \Pa{\sqrt{H}} + E \Pa{\sqrt{H}}}, \\
  \der A_2 &= \frac{\sqrt2}{\pi} K \Pa{\sqrt{H}} \, \der J + \frac{\sqrt2}{\pi} E \Pa{\sqrt{H}} \, \der H.
\end{align*}
In particular, $A_2$ is smooth there.
We also find the limit
\begin{align*}
  \lim_{(J, H) \to (0, 1)} A_2 &= \frac{2 \sqrt2}{\pi},
\end{align*}
which is the ``height'' of the $\times$ in \cref{fig:action-coord}.

\subsection{Action coordinates of a general spherical pendulum}

Now we come back to the general case where the potential is merely a function $V$ with maximum $V(-1) = V(1) = 1$ and minimum $V(0) = 0$, strictly decreasing in $[-1, 0]$ and strictly increasing in $[0, 1]$.
The procedures we take to get action coordinates in this case are analogous to the ones we take in the $V (z) = z^2$ case.
Using \cref{eq:X-A2-in-JH}, we get 
\begin{align*}
  \der A_1 &= \der J, \\
  \der A_2 &= \frac{\der J}{\pi}  \int_{z^-_{J, H}}^{z^+_{J, H}} \frac{J \der z}{(1 - z^2) \sqrt{2 (1 - z^2) (H - V (z)) - J^2}} - \frac{\der H}{\pi} \int_{z^-_{J, H}}^{z^+_{J, H}} \frac{\der z}{\sqrt{2 (1 - z^2) (H - V (z)) - J^2}},
\end{align*}
where $z = z^\pm_{J, H}$ are the unique nonnegative/nonpositive solutions to 
\begin{equation*}
  (1 - z^2) (H - V (z)) = \frac{J^2}{2}.
\end{equation*} 
Integrating $\der A_1$ and $\der A_2$, we obtain
\begin{equation} \begin{split} \label{eq:action-spdlm-gen}
  A_1 &= J, \\
  A_2 &= -\frac{1}{\pi} \int_{z^-_{J, H}}^{z^+_{J, H}} \frac{\sqrt{2(1 - z^2) (H - V (z)) - J^2}}{1 - z^2} \, \der z.
\end{split} \end{equation} 

\appendix

\section{Elliptic integrals} \label{sec:elliptic-integral}

In this appendix, we give a short introduction to the elliptic integrals used in our computation.
 
The \emph{elliptic integral of the first kind} is defined as
\begin{equation} \label{eq:ell-int-1}
  F (\gamma, k) = \int_0^\gamma \frac{1}{\sqrt{1 - k^2 \sin^2 \gamma}} \der \gamma,
\end{equation}
and the \emph{complete elliptic integral of the first kind} may thus be defined as 
\begin{equation} \label{eq:ell-int-1-comp}
  K (k) = F \Pa{\frac\pi2, k} = \int_0^{\frac\pi2}\frac{\der \gamma}{\sqrt{1 - k^2 \sin^2 \gamma}}.
\end{equation}
The values of $K$ at $k = 0$ and the asymptotics of $K$ at $k = 1$ are
\begin{equation} \begin{split} \label{eq:ell-int-1-comp-asym}
  K (0) &= \frac{\pi}{2}, \\
  K (k) &= -\frac12 \ln (1 - k) + \frac32 \ln 2 - \frac14 (1 - k) \ln (1 - k) + \mathcal{O} (1 - k) \text{ as } k \to 1-.
\end{split} \end{equation}

The \emph{elliptic integral of the second kind} is defined as
\begin{equation} \label{eq:ell-int-2}
  E (\gamma, k) = \int_0^\gamma \sqrt{1 - k^2 \sin^2 \gamma} \der \gamma,
\end{equation}
and the \emph{complete elliptic integral of the second kind} may thus be defined as 
\begin{equation} \label{eq:ell-int-2-comp}
 E (k) = E \Pa{\frac\pi2, k} = \int_0^{\frac\pi2} \sqrt{1 - k^2 \sin^2 \gamma} \der \gamma,
\end{equation}
with values at $k = 0$ and $k = 1$ as
\begin{equation} \label{eq:ell-int-2-comp-asym}
  E (0) = \frac{\pi}{2}, \qquad E (1) = 1.
\end{equation}

The \emph{elliptic integral of the third kind} is defined as
\begin{equation} \label{eq:ell-int-3}
  \Pi (\gamma, n, k) = \int_0^\gamma \frac {1}{(1 - n \sin^2 \gamma) \sqrt{1 - k^2 \sin^2 \gamma}} \der \gamma,
\end{equation}
and the \emph{complete elliptic integral of the third kind} is
\begin{equation} \label{eq:ell-int-3-comp}
  \Pi (n, k) = \Pi \Pa{\frac{\pi}{2}, n, k} = \int_0^{\frac\pi2}\frac{\der \gamma}{\Pa{1 - n \sin^2 \gamma}\sqrt{1 - k^2 \sin^2 \gamma}}.
\end{equation}
The asymptotics of $\Pi$ at $n = 0$ and $n = 1$ are
\begin{equation} \begin{split} \label{eq:ell-int-3-comp-asym}
  \Pi (n, k) &= K (k) - \frac{E (k) - K (k)}{k^2} n + \mathcal{O} (n^2) \text{ as } n \to 0+, \\
  \Pi (n, k) &= \frac{\pi}{2 \sqrt{(1 - \frac{k^2}{n}) (1 - n)}} + \frac{E (k)}{k^2 - 1} + K (k) + \mathcal{O} (1 - n) \text{ as } n \to 1-.
\end{split} \end{equation}

The \emph{Jacobi amplitude} is the inverse of the elliptic integral of the first kind:
\begin{equation} \label{eq:jacobi-am}
  f = F (\operatorname{am}(f, k), k),
\end{equation}
and the \emph{Jacobi elliptic sine} is the sine of the Jacobi amplitude given by
\begin{equation} \label{eq:jacobi-sn}
  \operatorname{sn}(f, k) = \sin \operatorname{am}(f, k).
\end{equation}

\end{document}